\documentclass[aap,MSNbibl,citesort,dvips]{arximspdf}
\usepackage{mathbh}

\doi{10.1214/10-AAP751}
\volume{22}
\issue{2}
\pubyear{2012}
\firstpage{457}
\lastpage{476}

\makeatletter

\newtheorem{Lemma}{Lemma}
\newtheorem{Proposition}[Lemma]{Proposition}
\newtheorem{Theorem}[Lemma]{Theorem}

\newproclaim{rema}{Remark}[section]
\newproclaim{rem}{Preliminary Remark}

\newcommand{\implies}{\Longrightarrow}
\newcommand{\AAA}{\mathcal{A}}
\newcommand{\des}{\mathop{\longrightarrow}\limits}
\newcommand{\one}{\mathbh1}
\newcommand{\iy}{\infty}
\newcommand{\R}{\mathbb{R}}
\newcommand{\N}{\mathbb{N}}
\newcommand{\Cc}{\mathcal{C}}
\newcommand{\F}{\mathcal{F}}
\newcommand{\Es}{\mathbb{E}}
\newcommand{\Pb}{\mathbb{P}}
\newcommand{\Eb}{\mathbb{E}}
\newcommand{\Ff}{\mathbb{F}}
\newcommand{\C}{\mathcal{C}}
\newcommand{\A}{\mathcal{A}}

\newcommand{\al}{\alpha}
\newcommand{\be}{\beta}
\newcommand{\G}{\Gamma}
\newcommand{\De}{\Delta}
\newcommand{\e}{\varepsilon}
\newcommand{\et}{\eta}
\newcommand{\tet}{\theta}
\newcommand{\La}{\Lambda}

\makeatother

\setattribute{abstract}{width}{295pt}

\begin{document}
\begin{frontmatter}

\title{On ergodic two-armed bandits}
\runtitle{Ergodic two-armed bandits}

\begin{aug}
\author[A]{\fnms{Pierre} \snm{Tarr\`{e}s}\corref{}\thanksref{t1}\ead[label=e1]{tarres@math.univ-toulouse.fr}\ead[label=u1,url]{http://www.math.univ-toulouse.fr/\textasciitilde tarres/}} and
\author[B]{\fnms{Pierre} \snm{Vandekerkhove}\ead[label=e2]{pierre.vandek@univ-mlv.fr}}
\runauthor{P. Tarr\`{e}s and P. Vandekerkhove}
\affiliation{CNRS, Universit\'{e} de Toulouse and Universit\'{e} Paris-Est}
\address[A]{Institut de Math\'{e}matiques \\
CNRS, Universit\'{e} de Toulouse\\
118 route de Narbonne\\
31062 Toulouse Cedex 9\\
France\\
\printead{e1}} %adresu isvedimo komanda gale!
\address[B]{Universit\'{e} Paris-Est\\
Marne-la-Vall\'{e}e, LAMA\\
5 boulevard Descartes \\
Champs-sur-Marne\\
77454 Marne-la-Vall\'{e}e Cedex 2\\
France\\
\printead{e2}}
\end{aug}

\thankstext{t1}{On leave from the University of Oxford.
Supported in part by the Swiss National Foundation Grant
200021-1036251/1 and by a Leverhulme Prize.}

% HISTORY:
\received{\smonth{5} \syear{2009}}
\revised{\smonth{11} \syear{2010}}

% ABSTRACT
%
\begin{abstract}
A device has two arms with \textit{unknown deterministic payoffs} and the
aim is to asymptotically identify the best one without spending too
much time on the other. The Narendra algorithm offers a stochastic
procedure to this end. We show under weak ergodic assumptions on these
deterministic payoffs that the procedure eventually chooses the best
arm (i.e., with greatest Cesaro limit) with probability one for
appropriate step sequences of the algorithm. In the case of i.i.d.
payoffs, this implies a ``quenched'' version of the ``annealed'' result
of Lamberton, Pag\`{e}s and Tarr\`{e}s [\textit{Ann. Appl. Probab.}
\textbf{14} (2004) 1424--1454] by the law of iterated logarithm, thus
generalizing it.

More precisely, if $(\eta_{\ell,i})_{i\in\N}\in\{0,1\}^\N$,
$\ell\in\{A,B\}$, are the deterministic reward sequences we would get if
we played at time $i$, we obtain infallibility with the same assumption
on nonincreasing step sequences on the payoffs as in Lamberton,
Pag\`{e}s and Tarr\`{e}s [\textit{Ann. Appl. Probab.} \textbf{14}
(2004) 1424--1454], replacing the i.i.d. assumption by the hypothesis
that the empirical averages $\sum_{i=1}^n\eta_{A,i}/n$ and
$\sum_{i=1}^n\eta_{B,i}/n$ converge, as $n$ tends to infinity,
respectively, to $\theta_A$ and $\theta_B$, with rate at least $1/(\log
n)^{1+\e}$, for some $\e>0$.

We also show a fallibility result, that is, convergence with positive
probability to the choice of the wrong arm, which implies the
corresponding result of Lamberton, Pag\`{e}s and Tarr\`{e}s
[\textit{Ann. Appl. Probab.} \textbf{14} (2004) 1424--1454] in the
i.i.d. case.
\end{abstract}

% KEYWORDS
%

\end{frontmatter}

%s1 ###
\section{Introduction}
\label{intro}
%s1.1 ###
\subsection{General introduction}
The so-called two-armed bandit is a device with two arms, each one
yielding an outcome in $\{0,1\}$ at each time step, irrespective of the
strategy of the player, who faces the challenge of choosing the best
one without losing too much time on the other.

The Narendra algorithm is a stochastic procedure devised to this end
which was initially introduced by Norman \cite{N68} and Shapiro and
Narendra~\cite{SN69} (see also \mbox{\cite{NT74,NT89}}) in the fields of mathematical psychology and
learning automata.
An application to optimal adaptive asset allocation in a financial
context has been developed by Niang \cite{N99}.

Formally, let $(\Omega, {\mathcal F}, \Pb)$ be a probability space.
The Narendra two-armed bandit algorithm is defined as follows. At each
time step $n\in\N$, we play source $A$ (resp., source $B$) with
probability $X_n$ (resp., $1-X_n$), where $X_0=x\in(0,1)$ is fixed and
$X_n$ is updated according to the following rule, for all $n\geq0$:
%
%e1 ###
\begin{equation}\label{bandit}
X_{n+1}=\cases{
X_n+\gamma_{n+1}(1-X_n), &\quad if $U_{n+1}=A$ and $\eta_{A,n+1}=1$,\cr
(1-\gamma_{n+1})X_n, &\quad if $U_{n+1}=B$ and $\eta_{B,n+1}=1$,\cr
X_n, &\quad otherwise,}
\end{equation}
where $(\gamma_n)_{n\geq1}$ is a deterministic sequence taking values
in $(0,1)$, $U_{n+1}$ is the random variable corresponding to the
label of the arm played at time $n+1$ and $\eta_{\ell,n+1}$ denotes
the payoff, taking values in $\{0,1\}$, of source $\ell
\in\{A,B\}$ at time $n+1$.

We assume without loss of generality that $U_{n+1}=A \one_{\{
I_{n+1}\leq
X_n\}}+\break B\one_{\{I_{n+1}>X_n\}}$, where $(I_n)_{n\geq1}$ is a sequence of
independent uniformly distributed random variables on $[0,1]$.

The literature on this algorithm generally assumes that the sequences
$(\eta_{A,n})_{n\geq1}$ and $(\eta_{B,n})_{n\geq1}$ are independent
with Bernoulli distributions of parameters $\tet_A$ and~$\tet_B$,
where $\tet_A>\tet_B$, the aim being to determine whether
$(X_n)_{n\in\N}$ a.s. converges to~$1$ or not as $n$ tends to infinity.

Notwithstanding the apparent simplicity of this stochastic procedure,
the first criteria on a.s. convergence to ``the good arm'' under the
above i.i.d. assumptions were only obtained thirty years after the
original definition of this Narendra algorithm by Tarr\`{e}s
\cite{T01} and Lamberton, Pag\`{e}s and Tarr\`{e}s \cite{LPT04} in a more
general framework.
Recently Lamberton and Pag\`{e}s established the corresponding rate of
convergence \cite{LP05a} and proposed and studied a penalized version
\cite{LP05b}. Note that a game theoretical question arising in the
context of two-armed bandits was recently studied by Bena\"{\i}m and
Ben Arous \cite{BB03} and Pag\`{e}s \cite{P05}.

Our work focuses on the understanding of the Narendra two-armed bandit
algorithm under the assumption that the payoff sequences $(\eta_{\ell
,n})_{n\geq1}$, $\ell\in\{A,B\}$, are \textit{unknown} and
\textit{deterministic}. Under the following condition (S) on the
step sequence (required in \cite{LPT04} but without monotonicity) and
weak ergodic assumption (E2) on the rate at which $A$
must be asymptotically better than $B$, we show that $X_n$ a.s.
converges to $1$. Heuristically, the result points out that, even with
strongly dependent outcomes, $X_n$ accumulates sufficient statistical
information on the ergodic behavior of the two arms to induce a
corresponding appropriate decision.

More precisely, let us introduce the following \textit{step sequence} and
\textit{ergodic} assumptions.

\subsubsection*{Step sequence conditions}
Let, for all $n\in\N\cup\{\iy\}$,
$ \G_n=\sum_{k=1}^n \gamma_k$.

Let (S1) and (S2) be the following assumptions on the step
sequence $(\gamma_n)_{n\in\N}$:
\begin{longlist}[(S2)]
\item[(S1)] $(\gamma_n)_{n\geq1}$ is nonincreasing and $\G_\iy=\iy$;
\item[(S2)] $\gamma_n=O(\G_ne^{-\tet_B\G_n})$.
\end{longlist}
Let (S) be the set of conditions (S1) and (S2).

\subsubsection*{Ergodic conditions}
Let (E) be the assumption that the ouputs of arms~$A$ and $B$ satisfy
\[
\mbox{(E)\qquad\qquad\quad\hphantom{,,}}\frac{1}{n}\sum_{k=1}^n\eta_{A,k}\des
_{n\rightarrow\infty}\theta_A \quad\mbox{and}\quad \frac{1}{n}\sum
_{k=1}^n\eta_{B,k}\des_{n\rightarrow\infty}\theta_B,\qquad\qquad\quad
\]
where $\theta_A$, $\theta_B$ $\in(0,1)$.
The ergodic condition (E) means that the average payoff of arm
$A$ (resp., arm $B$) is $\tet_A$ (resp., $\tet_B$) but does not
assume anything on the corresponding rate of convergence. In order to
introduce conditions on this rate, let us denote, for all $n\in\N$,
\[
R_n:=\max_{\ell\in\{A,B\}}\Biggl|\sum_{i=1}^n (\eta
_{\ell,i}-\theta_\ell)\Biggr|.
\]

Given a map $\phi\dvtx\N\longrightarrow\R_+$ and $\theta_A$, $\theta
_B$ $\in(0,1)$, let us denote by

\begin{longlist}[(E$\phi$)]
\item[(E$\phi$)]
the assumption that $R_n/\phi(n) \des
_{n\rightarrow\infty}0$.
\end{longlist}

Let (E1) and (E2) be condition (E$\phi$),
respectively, with the following assumption on $\phi$:

\begin{longlist}[(E2)]
\item[(E1)]
$\phi$ is nondecreasing concave on $[k_0,\iy)$
for some $k_0\in\N$ and
\[
\sup_{n\in\N} \gamma_n \phi(n)<\infty.
\]
\item[(E2)]
$ \phi(n)=\frac{n}{(\log(n+2))^{1+\varepsilon}}$ for some $\varepsilon>0$.
\end{longlist}

Note that (E) corresponds to (E$\phi$) with $\phi(n)=n$,
$n\in\N$, under which~(E1) holds, for instance, in the case of a
step sequence $\gamma_n=c/(c+n)$, $c>0$. Also, Lemma \ref{ubound},
proved in Section \ref{detest}, implies that (S)--(E2)
$\implies$ (E1).
\begin{Lemma}\label{ubound} If condition \textup{(S)} holds, then
\[
\limsup_{n\to\iy}\frac{\gamma_nn}{\log
n}\leq\limsup_{n\to\iy}\frac{\G_n}{\log n}\leq1/\tet_B.
\]
\end{Lemma}

Theorems \ref{thm:conv} and \ref{thm:rightconv} provide assumptions
for convergence of the Narendra sequence $(X_n)_{n\geq0}$ toward $0$ or
$1$ as $n$ tends to infinity, respectively, convergence toward $1$
when $\tet_A>\tet_B$ (i.e., asymptotic choice of the ``right arm'').

\begin{Theorem}\label{thm:conv} Under assumptions \textup{(S1)--(E1)},
the Narendra
sequen-\break ce~$(X_n)_{n\in
\mathbb{N}}$ converges $\mathbb{P}_x$-a.s. toward $0$ or $1$ as $n$
tends to infinity.
\end{Theorem}
\begin{Theorem}\label{thm:rightconv} Under assumptions \textup{(S)--(E2)}
and $\tet
_A>\tet_B$, the Narendra sequence
$(X_n)_{n\in\mathbb{N}}$ converges $\mathbb{P}_x$-a.s. toward $1$
as $n$ tends to infinity.
\end{Theorem}

Recall that the above conditions (E1) and (E2) are purely
deterministic. If we let the sequences
$(\eta_{A,i})_{i\in\N}$ and $(\eta_{B,i})_{i\in\N}$ be
distributed as i.i.d. sequences with expectations
$\tet_A$ and $\tet_B$, then (E2) almost surely occurs as
a~consequence of the law of iterated logarithm.
Assuming (S) and $\tet_A>\tet_B$, Theorem \ref{thm:rightconv}
implies that the algorithm $(X_n)_{n\in\N}$ almost surely
converges to~$1$, which is a generalization of the
corresponding infallibility Proposition 5 proved by Lamberton, Pag\`
{e}s and Tarr\`{e}s in \cite{LPT04} for nonincreasing step sequences
$(\gamma_n)_{n\in\N}$.

In practice, the Narendra algorithm is used in the context of
performance assessment, or in applications either in automatic control or
in financial mathematics and the i.i.d. assumption looks rather
unrealistic since the performance depends in general on parameters that
evolve slowly and randomly in time. The following framework provides a
possible generalization.

Suppose that $(S_{\ell,i}) _{i\in\N}$, $\ell\in\{A,B\}
$, are ergodic stationary Markov chains taking values in a measurable
space $(\mathbb X, \mathcal X)$, with transition kernel $Q_\ell$ and
stationary initial distribution $\pi_\ell$. Let us consider a
measurable event $\Cc\in\mathcal X$, and define sequences
$(\eta_{\ell,i})_{i\in\N}$, for $\ell\in\{A,B\}$, as
%
%e2 ###
\begin{equation}
\label{mtarget}
\eta_{\ell,i}=\one_{\{S_{\ell,i}\in\Cc\}}, \qquad  i \in\N.
\end{equation}
These random sequences $(\eta_{\ell,i})_{i\in\N}$ are functions of the
states of the Markov chains and satisfy, as a consequence, the ergodic
condition (E), with
\[
\theta_\ell=\pi_\ell(S_{\ell,0}\in\Cc).
\]

The sequences $(S_{\ell,i}) _{i\in\N}$, $\ell\in\{A,B\}$,
represent the agents' outputs from which $(\eta_{\ell,i})_{i\in\N}$
extracts scores through target assessment. Note that, contrary to
$(S_{\ell,i}) _{i\in\N}$, $(\eta_{\ell,i})_{i\in\N}$ is not
Markov in general.

Miao and Yang \cite{MY08} establish under weak conditions (concerning
mainly the transition kernels $Q_\ell$) the law of iterated logarithm
for additive functionals of Markov chains, thus providing the required
ergodic rate of convergence~(E2).

Let us now show a simple fallibility result that will also imply the
corresponding result of \cite{LPT04} in the i.i.d. case.
\begin{Theorem}
\label{fallibility}
Assume $\tet_A>\tet_B$ and $\sum_{n\geq0}\prod_{k=1}^n(1-\gamma
_k\eta_{B,k})<\iy$. Then $\Pb(\lim_{n\to\iy}X_n=0)>0$.\vadjust{\goodbreak}
\end{Theorem}
\begin{rema}
In the case where $(\eta_{B,k})_{k\geq0}$ is an i.i.d. sequence of
random variables, then
\[
\Es_x\Biggl(\sum_{n\geq0}\prod_{k=1}^n(1-\gamma_k\eta_{B,k})\Biggr)=
\sum_{n\geq0}\prod_{k=1}^n(1-\gamma_k\tet_B)<\iy
\]
ensures that the third condition of Theorem \ref{fallibility} is
fulfilled and, therefore, Theorem~\ref{fallibility} implies the
fallibility result Theorem 1(b) in \cite{LPT04}.
\end{rema}
\begin{rema}
In the general (ergodic) case, if $\sum\gamma_n^2<\iy$, $\sum\G
_n|\phi''(n)|<\iy$ and $\limsup\G_n|\phi'(n)|<\iy$, then the
proof of Lemma \ref{Abel2} implies that the conditions of Theorem
\ref{fallibility} are equivalent to $\sum\exp(-\G_n\tet_B)<\iy$ and
$\tet_A>\tet_B$. These assumptions hold, for instance, if $\gamma
_n=c/(c+n)$ and $\phi(n)=n/(\log(n+2))^{1+\e}$ for some $\e>0$ and
$c\tet_B>1$ (see also the proof of Lemma \ref{Abel2}).
\end{rema}
\begin{pf*}{Proof of Theorem \ref{fallibility}}
Recall that $X_0=x\in(0,1)$. Let $A$ be the event
\[
A:=\{\forall k\geq1,  I_k\leq X_k\}
=\Biggl\{\forall n\geq0,  X_n=x\prod_{k=1}^n(1-\gamma_k\eta
_{B,k})\Biggr\}.
\]
Then
\[
\Pb(A)=\prod_{n=1}^\iy\Biggl(1-x\prod_{k=1}^n(1-\gamma_k\eta
_{B,k})\Biggr)>0
\quad\iff\quad\sum_{n\geq0}\prod_{k=1}^n(1-\gamma_k\eta_{B,k})<\iy,
\]
and note that this last predicate, which is the second assumption of
the theorem, obviously implies $\sum\gamma_n\eta_{B,n}=\iy$.
Now, a.s. on $A$,
\[
X_n\leq x\exp\Biggl(-\sum_{k=1}^n\gamma_n\eta_{B,n}\Biggr)\des
_{n\to\iy}0,
\]
which concludes the proof.
\end{pf*}

\subsubsection*{Notation}
The letter $C$ will denote a positive real constant that may change
from one inequality to the other.

We write $\phi'$ and $\phi^{\prime\prime}$ for the first- and second-order
discrete derivatives of~$\phi$: for all $n\geq1$,
\[
\phi'(n):=\phi(n)-\phi(n-1)\quad  \mbox{and}\quad
\phi''(n):=\phi(n-1)+\phi(n+1)-2\phi(n).
\]

We let, for all $n\in\N$,
\[
\alpha_n:=R_n/\phi(n),\qquad \beta_n:=\sup_{k\geq n}\alpha_k.
\]
Note that, under assumption (E$\phi$), $\al_n$, $\be_n$ $\des
_{n\rightarrow\infty}0$.

Given two real sequences $(u_n)_{n\geq0}$ and $(v_n)_{n\geq0}$, we
write
\[
u_n=\Box(v_n),
\]
when, for all $n\geq0$, $|u_n|\leq|v_n|$.

%s1.2 ###
\subsection{\texorpdfstring{Sketch of the proofs of Theorems \protect\ref{thm:conv} and \protect\ref{thm:rightconv}}
{Sketch of the proofs of Theorems 2 and 3}}
\label{sec:sketch}

Our first aim is to write down in Proposition \ref{evol} the evolution
of $(X_n)_{n\geq0}$ as a stochastic perturbation of the Cauchy--Euler
procedure defined by
%
%e3 ###
\begin{equation}
\label{eq:cauchy}
x_{n+1}=x_n+\gamma_{n+1}h(x_n),
\end{equation}
where $h(x):=(\tet_A-\tet_B)f(x)$, with $f(x):=x(1-x)$.

However, contrary to the case of i.i.d. payoff sequences $(\eta_{\ell
,n})_{n\geq0}$, $\ell\in\{A,B\}$, considered in \cite{LPT04}, the
perturbation of the scheme (\ref{eq:cauchy}) under an ergodic
assumption (E) does not only consist of a martingale, but also of
an increment whose importance depends on $\phi$, that is, on the rate
of convergence of the mean payoffs to $\tet_A$ and $\tet_B$. More
precisely let, for all $n\geq1$,
\[
\wedge_{n}=\sum_{k=1}^n \gamma_k f(X_{k-1})\bigl(\eta_{A,k}-\eta
_{B,k}-(\theta_A-\theta_B)\bigr)
\]
with the convention that
$\wedge_0=0$
and let $(M_n)_{n\geq1}$ be an $(\F_n)_{n\geq1}$-adapted martingale
given by
\[
M_n:=\sum_{k=1}^n\gamma_k\e_k,\qquad M_0:=0
\]
with
\[
\e_k:=
\eta_{A,k}(1-X_{k-1})(\one_{U_{k}=A}-X_{k-1})+\eta
_{B,k}X_{k-1}\bigl((1-X_{k-1})-\one_{U_{k}=B}\bigr).
\]

\begin{Proposition}
\label{evol}
For all $n\in\N$,
\[
X_n=x+M_n+\wedge_n+(\tet_A-\tet_B)\sum_{k=1}^n\gamma_kf(X_{k-1}).
\]
\end{Proposition}
\begin{pf}
The updating rule (\ref{bandit}) can be rewritten as
%
%e4 ###
\begin{eqnarray}\label{bandit_bis}
X_{n+1}&=&X_n+\gamma_{n+1} \eta_{A,n+1}(1-X_n)\one
_{U_{n+1}=A}-\gamma_{n+1}\eta_{B,n+1}X_n\one_{U_{n+1}=B}\nonumber\\[-2pt]
&=&X_n+\gamma_{n+1}
\eta_{A,n+1}(1-X_n)(\one_{U_{n+1}=A}-X_n)\nonumber\\[-9.5pt]\\[-9.5pt]
&&{} +\gamma_{n+1} \eta_{B,n+1}X_n\bigl((1-X_n)-\one_{U_{n+1}=B}\bigr)\nonumber\\[-2pt]
&&{} +\gamma
_{n+1} f(X_n)(\eta_{A,n+1}-\eta_{B,n+1}).\nonumber
\end{eqnarray}
\upqed\end{pf}

Note that Proposition \ref{evol} can be interpreted as the property
that the noise is multiplicative in the sense that, for all $n$,
\[
\gamma_{n+1}^{-1}(\La_{n+1}-\La_n)=f(X_n)\bigl(\eta_{A,k}-\eta
_{B,k}-(\theta_A-\theta_B)\bigr)
\]
is the product of a function of $X_n$ and a function of $(\eta
_{A,n+1},\eta_{B,n+1})$ outcome of the two arms at time $n+1$.

Let us now provide estimates of the evolution of $(\wedge_n)_{n\in\N
}$, which will be necessary to the proof of Theorem \ref
{thm:rightconv}; they will also imply Theorem \ref{thm:conv} in
passing. We note that Laruelle and Pag\`{e}s \cite{LP10} recently
generalized the proof\vadjust{\goodbreak} of this latter result as convergence of the
ergodic dynamics toward an equilibrium point of the corresponding ODE
under the assumption that the noise is multiplicative and a classical
Lyapounov assumption, or more generally under a strong Lyapounov
assumption, and technical conditions.

%Now, $\eta_{A,k}-\eta_{B,k}-(\theta_A-\theta_B)$ being erratic by
%the very nature of the question, these bounds will naturally involve a
%discrete integration by parts in order to make the ergodic upper bound
%function $\phi$ appear. However, $(\gamma_nf(X_{n-1}))_{n\in\N}$
%is not a nonincreasing sequence in general so that the technique cannot
%work directly.

Our estimates of $\La_n-\La_m$ for large $m$ and $n$ are derived by
discrete integration by parts. To this end, we need to get round the
difficulty that the sequence $(\gamma_nf(X_{n-1}))_{n\in\N}$ is not
monotonic in general.

Instead, let us define, for all $n\in\N$,
\[
\Delta_n:=\frac{\gamma_n}{\prod_{k=1}^n(1-\gamma_k)},\qquad
S_n:=\frac{1}{\prod_{k=1}^n(1-\gamma_k)}
\]
with the convention that $\De_0=S_0:=1$. Remark that $S_n\rightarrow
\infty$ if and only if $\sum_{ n\geq1 } \gamma_n=+\infty$.\vspace*{1pt}

Note that $x/S_n$ is a trivial lower bound for $X_n$ and that
%
%e5 ###
\begin{equation}\label{Delta-gamma}
\gamma_n=\frac{\Delta_n}{S_n}\qquad
\mbox{with }  S_n=\sum_{k=0}^n\Delta_k.
\end{equation}
We first study the sequence $(\Psi_n)_{n\in\N}$ defined by
\[
\Psi_n:=\sum_{k=n+1}^\iy\frac{\gamma_k}{S_{k-1}}\bigl(\et_{A,k}-\et
_{B,k}-(\tet_A-\tet_B)\bigr);
\]
$(\Psi_n)_{n\geq1}$ is well defined since, for all $\ell\in\{
A,B\}$,
\[
\sum_{k=2}^\iy\frac{\gamma_k}{S_{k-1}}|\et_{\ell,k}-\tet_\ell|
\leq\sum_{k=2}^\iy\frac{\gamma_k}{S_{k-1}}
=\sum_{k=2}^\iy\biggl(\frac{1}{S_{k-1}}- \frac{1}{S_{k}}
\biggr)=\frac{1}{S_1}
\]
since under (S1) we have $S_n\des_{n\rightarrow\infty} \infty$.
Since $(\gamma_n/S_{n-1})_{n\in\N}$ is a nonincreasing sequence if
$(\gamma_n)_{n\in\N}$ is itself nonincreasing [recall that $\gamma
_n\in(0,1)$], we deduce Lemma \ref{1sumest} by an Abel transform,
that is, discrete integration. Moreover, we observe that, for all $n\geq
m\geq0$, the evolution of $\wedge_{\cdot}$ between time steps $m$ and $n$ is
given by
\[
\wedge_{n}-\wedge_{m}=
\sum_{k=m+1}^n S_{k-1}f(X_{k-1})\frac{\gamma_k}{S_{k-1}}\bigl(\eta
_{A,k}-\eta_{B,k}-(\theta_A-\theta_B)\bigr).
\]

Now, $(S_kf(X_k))_{k\in\N}$ is a nondecreasing sequence. Indeed, for
all $k\in\N$, $f(X_k)\geq(1-\gamma_k)f(X_{k-1})$ since $f$ is
concave and $X_k$ is the barycentre of~$X_{k-1}$ and either $0$ or $1$,
with weights $1-\gamma_k$ and $\gamma_k$, where $f(0)=f(1)=0$. We
rely on this monotonicity and apply an Abel transform again, which
enables us to show Lemma \ref{2sumest}.
\begin{Lemma}
\label{1sumest}
Assume that $(\gamma_n)_{n\in\N}$ is nonincreasing and that $\phi$
is nondecreasing concave on $[k_0,\iy)$ for some $k_0\in\N$. Then,
for all $n\geq k_0$,
\[
|\Psi_n|\leq\frac{2\beta_n}{S_{n-1}}[\phi'(n)+2\gamma_n
\phi(n)].\vadjust{\goodbreak}
\]
\end{Lemma}
\begin{Lemma}
\label{2sumest}
Let, for all $n\in\N$,
\[
R'_n:=\frac{2\sup_{k\geq n}\beta_k[\phi'(k)+2\gamma_k \phi(k)
]}{1-\gamma_n}.
\]
Under the assumptions of Lemma \ref{1sumest} we have, for all $n\geq
m\geq k_0$,
\[
|{\wedge_{n}-\wedge_{m}}|\leq R'_m\Biggl[\sum_{k=m+1}^n\gamma
_kf(X_{k-1})+2f(X_n)\Biggr].
\]
\end{Lemma}

Lemmas \ref{1sumest} and \ref{2sumest} are proved in Sections \ref
{sec:1sumest} and \ref{sec:2sumest}.

These results enable us to conclude the proof of Theorem \ref
{thm:conv}. Indeed, by Proposition \ref{evol} and Lemma \ref
{2sumest}, for all $n\geq m\geq0$,
%
%e6 ###
\begin{eqnarray}
\label{evolm}
X_n-X_m&=&M_n-M_m+\wedge_n-\wedge_m+(\tet_A-\tet_B)\sum_{k=m+1}^n
\gamma_kf(X_{k-1})\nonumber\\
&=&M_n-M_m+\bigl(\tet_A-\tet_B+\Box(R'_m)\bigr)\sum_{k=m+1}^n \gamma
_kf(X_{k-1})\\
&&{}+2\Box(R'_m) f(X_n).\nonumber
\end{eqnarray}
We assume that (E1) and (S1) hold; thus, $R'_n\des_{n\to
\iy}0$.
Let us prove by contradiction that
%
%e7 ###
\begin{equation}
\label{sumfin}
\sum_{k=1}^\infty\gamma_k f(X_{k-1})<\infty\qquad\mbox{a.s.}
\end{equation}
holds. Indeed, let us assume the contrary; choose $m$ such that
$|R'_m|<|\tet_A-\tet_B|$. A.s.\vspace*{1pt} on
$\{\sum_{k=1}^\infty\gamma_k f(X_{k-1})=\infty\}$, using Chow's
lemma (see, e.g., \cite{Duflo}) and $\Es(\e_{k+1}^2|\F
_k)\leq2 f(X_k)$,
we deduce
\[
M_n-M_m=o\Biggl(\sum_{k=m+1}^n\gamma_k^2f(X_{k-1})\Biggr)=o
\Biggl(\sum_{k=m+1}^n\gamma_kf(X_{k-1})\Biggr)
\]
and, therefore, for all $n$, $m$ $\in\N$,
\[
X_n-X_m=\bigl(\tet_A-\tet_B+\Box(R'_m)+o_{n\to\iy}(1)\bigr)\sum
_{k=m+1}^n\gamma_kf(X_{k-1})+O(1),
\]
which is contradictory using $X_n\in[0,1]$ for all $n\in\mathbb
N$.

Hence, $\mathbb P_x$-almost surely,
$(X_n)_{n\geq0}$ is a Cauchy sequence and, therefore, converges to a
limit random variable $X_\iy\in[0,1]$. Now (\ref{sumfin}) implies
that $f(X_\iy)=0$, since $\G_\iy=\iy$ and, therefore, $X_\iy=0$ or
$1$ a.s.

The proof of Theorem \ref{thm:rightconv} itself has two parts. The
first one consists in showing a ``brake phenomenon,'' that is, that
$(X_n)_{n\geq0}$ cannot in any case decrease too rapidly to $0$ as $n$
goes to infinity. We already observed that, trivially, $X_n$ is lower
bounded\vadjust{\goodbreak} by $x/S_n$. A better lower bound can easily be obtained; let us
define, for all $n\in\N$,
\[
S_n^B:=\frac{1}{\prod_{k=1}^n(1-\gamma_k
\one_{\{I_k>X_{k-1},\eta_{B,k}=1\}})}\qquad \mbox{with initial condition }
S_0^B=0
\]
and, for all $n\geq1$,
\[
\De_n^B:=\gamma_nS_n^B,\qquad Y_n^B:=S_n^BX_n.
\]

Note that, as a consequence of the definition of the Narendra
algorithm~(\ref{bandit}), for all $n\geq0$,
%
%e8 ###
\begin{equation}\label{ynbupdate}
Y_{n+1}^B=\cases{
Y_n^B+\De_{n+1}^B(1-X_n), &\quad if $U_{n+1}=A$ and $\eta_{A,n+1}=1$,\cr
Y_n^B, &\quad otherwise.}
\end{equation}

Roughly\vspace*{1pt} speaking, $S_n^B$ is the product $S_n$ restricted to playing
and winning with $B$; $x/S_n^B$ is straightforwardly a lower bound of
$X_n$. Proposition \ref{pa=0}, proved in Section \ref{sec:pa=0},
further claims that, for any $C>0$, $C\log S_n^B/S_n^B$ is an
asymptotic lower bound of $X_n$ a.s. on $\{X_\iy=0\}$.\vspace*{-3pt}
\begin{Proposition}
\label{pa=0}
Under assumptions \textup{(S)} and \textup{(E2)},
\[
\Bigl\{\lim_{n\to\iy}X_n=0\Bigr\}\subseteq\biggl\{\limsup_{n\to
\iy}\frac{X_n}{\log S_n^B/S_n^B}=\iy\biggr\},\qquad \Pb_x\mbox{-a.s.}
\]
\end{Proposition}

The second part of the proof of Theorem \ref{thm:rightconv} assumes
$\tet_A>\tet_B$ and is given in Section \ref{sec:conc}. Recall that,
by Theorem \ref{thm:conv}, $X_n$ converges a.s. to $0$ or $1$ [using
the remark that (S)--(E2) implies (E1), see the remark before
the statement of Lemma~\ref{ubound}] so that we only need to show that
$\Pb(\lim X_n=0)=0$.

We study $(X_n)_{n\geq0}$ as a perturbed Cauchy--Euler scheme and prove
by Doob's inequality that, starting from $C\log S_n^B/S_n^B$ for
sufficiently large \mbox{$C>0$}, $X_n$ remains bounded away from $0$ with
lower bounded probability, which enables us to conclude that $X_\iy
\not=0$ a.s.\vspace*{-3pt}

%s2 ###
\section{Deterministic estimates on the step sequence}
\label{detest}
We first recall below the two following preliminary remarks in \cite
{LPT04} that (S2) implies on one hand that $\sum_{n=1}^\infty
\gamma_n^2<\iy$ and, on the other hand, that $\G_n-\log S_n$
converges as $n$ goes to infinity.

Then we prove Lemma \ref{ubound} that (S) implies explicit
asymptotic upper bounds on $(\gamma_n)_{n\in\N}$ and $(\G_n)_{n\in
\N}$.\vspace*{-3pt}
\begin{rem}\label{prrem1}
Assumption (S2) implies
$\sum_{n=1}^\infty\gamma_n^2<\iy$ since, for all $n\in\N$,
\begin{eqnarray*}
\sum_{k=1}^n\gamma_k^2&\leq& C\sum_{k=1}^n(\G_k-\G_{k-1})\G
_ke^{-\tet_B\G_k}\\[-2pt]
&\leq& C\int_0^{\G_n}ue^{-\tet_Bu} \,du\leq C\int_0^{+\iy}ue^{-\tet
_Bu} \,du<\iy
\end{eqnarray*}
using that $u\mapsto ue^{-\tet_Bu}$ is nonincreasing for
$u>\tet_B^{-1}$.\vadjust{\goodbreak}
\end{rem}
\begin{rem}\label{prrem2}
The partial sums $S_n$ and $\Gamma_n$ satisfy for every \mbox{$n\geq1$},
%
%e9 ###
\begin{equation}\label{IneqSGamma}
\log S_n
-\sum_{k=1}^n\frac{\gamma^2_k}{1-\gamma_k}\leq\Gamma_n\leq\log S_n.
\end{equation}
This follows from the easy comparisons
\[
\Gamma_n =\sum_{k=1}^n
\frac{\Delta_k}{S_k}
\cases{
\displaystyle \leq  \int_1^{S_n}\frac{du}{u} =
\log S_n,\vspace*{3pt}\cr
\displaystyle  = \sum_{k=1}^n \frac{S_{k-1}}{S_k}
\int_{S_{k-1}}^{S_k} \frac{du}{S_{k-1}}\geq
\sum_{k=1}^n
(1-\gamma_k) \int_{S_{k-1}}^{S_k} \frac{du}{u},\vspace*{3pt}\cr
\displaystyle \geq \log S_n - \sum_{k=1}^n
\frac{\gamma^2_k}{1-\gamma_k}.}
\]
\end{rem}
\begin{pf*}{Proof of Lemma \ref{ubound}}
The first inequality is elementary, since $\Gamma_n\geq n\gamma_n$,
using that $(\gamma_n)_{n\geq1}$ is a nonincreasing sequence by
(S1). By assumption (S2), for some $C>0$, for all $n\in\N$,
\[
C\geq\frac{\gamma_n e^{\tet_B\G_n}}{\G_n}.
\]
Using that $u\mapsto e^{\theta_B u}/u$ is increasing on $[1/\theta
_B,\infty)$ we obtain that, for sufficiently large $n_0\in\N$,
%
%e10 ###
\begin{equation}
\label{eqn:est}
C(n-n_0)\geq\int_{\G_{n_0}}^{\G_n} \frac{e^{\tet_B x}}{x}
\,dx\mathop{\sim}_{n\to\iy}\frac{e^{\tet_B\G_n}}{\tet_B\G_n}.
\end{equation}
Trivially, $\log(e^{\tet_B\G_n}/\tet_B\G_n)\sim_{n\to\iy}\tet
_B\G_n$, so that (\ref{eqn:est}) proves the second inequality.
\end{pf*}

%s3 ###
\section{Abel transforms}
%s3.1 ###
\subsection{Preliminary estimates}
Lemmas \ref{Abel1} and \ref{Abel2} estimate the error in replacing
the payoffs $\eta_{\ell,k}$ by their ``average success rate'' $\tet
_\ell$ in a sum weighted by a decreasing sequence $(\xi_n)_{n\in\N
}$, by the use of Abel transforms,\vadjust{\goodbreak} that is, discrete integrations by
parts. More precisely let, for all $n\in\N$ and $\ell\in
\{A,B\}$,
\[
\Phi_{n,\xi}^\ell=\sum_{k=1}^n \xi_k(\et_{\ell,k}-\tet_\ell)
\]
be the corresponding deviation. Lemma \ref{Abel1} upper bounds $|\Phi
_{n,\xi}^\ell-\Phi_{m,\xi}^\ell|$ for all \mbox{$n\geq m$}, whereas Lemma
\ref{Abel2} shows that $\Phi_{n,\xi}$ converges to a finite
value
under certain assumptions, which are fulfilled, for instance, when $\xi
:=\gamma$ and~(S)--(E2) hold.\vadjust{\goodbreak}

Lemma \ref{Abel1} is the main tool in the proof of Lemmas \ref
{1sumest} and \ref{2sumest} and the second part of Lemma \ref{Abel2}
will be useful in the proof of Proposition \ref{pa=0} providing
``brake phenonemon'' bounds.\vspace*{-3pt}
\begin{Lemma}\label{Abel1} Let $(\xi_n)_{n\in\N}$ be a
positive real-valued nonincreasing sequence. Assume $\phi$ is
nondecreasing on $[k_0,\iy)$ for some $k_0\in\N$, then, for all
$n\geq m\geq k_0$,
\[
|\Phi_{n,\xi}^\ell-\Phi_{m,\xi}^\ell|\leq\beta_m\Biggl(
\sum_{k=m+1}^{n} \xi_k\phi'(k)+2\xi_{m}\phi(m)\Biggr).
\]
\end{Lemma}
\begin{pf}Let, for all $n\in\N$ and $\ell\in
\{A,B\}$, $\kappa_n^\ell:= \sum_{k=1}^n
(\eta_{\ell,k}-\theta_\ell)$. If $n\geq m\geq k_0$,
then
%
%e11 ###
\begin{eqnarray}\label{eq1:abel1}
\Phi_{n,\xi}^\ell-\Phi_{m,\xi}^\ell&=& \sum_{k=m+1}^n
\xi_k(\eta_{\ell,k}-\theta_\ell)\nonumber\\[-2pt]
&=& \sum_{k=m+1}^n \xi_k(\kappa_k^\ell-\kappa_{k-1}^\ell)=
\sum_{k=m+1}^n \xi_k \kappa_k^\ell-\sum_{k=m}^{n-1}
\xi_{k+1}\kappa_k^\ell\\[-2pt]
&=& \sum_{k=m}^{n-1} (\xi_k-\xi_{k+1})
\kappa_k^\ell+\xi_n\kappa_n^\ell-\xi_{m}\kappa_m^\ell.\nonumber
\end{eqnarray}
Now, using that $(\xi_n)_{n\geq0}$ is nonincreasing,
%
%e12 ###
\begin{eqnarray}\label{eq2:abel1}\quad
&&\Biggl|\sum_{k=m}^{n-1}(\xi_k-\xi_{k+1})\kappa_k^\ell\Biggr|\nonumber\\[-2pt]
&&\qquad\leq\sum_{k=m}^{n-1}(\xi_k-\xi_{k+1}) R_k
=\sum_{k=m}^{n-1}(\xi_k-\xi_{k+1}) \alpha_k\phi(k)\nonumber\\[-9.5pt]\\[-9.5pt]
&&\qquad\leq\beta_m\sum_{k=m}^{n-1}(\xi_k-\xi_{k+1}) \phi(k)
=\beta_m\Biggl(\sum_{k=m}^{n-1} \xi_k\phi(k)-
\sum_{k=m+1}^{n} \xi_{k}\phi(k-1)\Biggr)\nonumber\\[-2pt]
&&\qquad=\beta_m\Biggl(\sum_{k=m+1}^{n} \xi_k
\bigl(\phi(k)-\phi(k-1)\bigr)+\xi_{m}\phi(m)-\xi_n\phi
(n)\Biggr).\nonumber
\end{eqnarray}
In summary, (\ref{eq1:abel1}) and (\ref{eq2:abel1}) imply
\begin{eqnarray*}
|\Phi_{n,\xi}^\ell-\Phi_{m,\xi}^\ell|&\leq&\beta
_m\Biggl(\sum_{k=m+1}^{n} \xi_k
\bigl(\phi(k)-\phi(k-1)\bigr)+2\xi_m\phi(m)\Biggr)\\[-2pt]
&=&\beta_m\Biggl(\sum_{k=m+1}^{n} \xi_k\phi'(k)+2\xi_{m}\phi
(m)\Biggr).
\end{eqnarray*}
\upqed\end{pf}
\begin{rema}
\label{ubdphip}
Under assumption (E2), that is, when $\phi(k):=k(\log
(k+2))^{-(1+\e)}$ for some $\e>0$, then
\[
\phi'(k)\leq\frac{1}{(\log(k+1))^{1+\varepsilon}}, \qquad  k\in
\N.
\]
Indeed, for all $x\in\R^+$,
\[
\biggl(\frac{d\phi}{dx}\biggr)(x)=\frac{1}{(\log
(x+2))^{1+\varepsilon}}-\frac{(1+\varepsilon)x}{(x+2)(\log
(x+2))^{2+\varepsilon}}
\]
and
\[
\phi'(k)\leq\sup_{x\in[k-1,k]}\biggl(\frac{d\phi}{dx}\biggr)(x).
\]
\end{rema}
\begin{Lemma}\label{Abel2} Given a
positive real-valued nondecreasing sequence $(\xi_n)_{n\in\N}$, let,
for all $n\in\N$, $\Xi_n:=\sum_{k=1}^n\xi_k$. If $\phi$ is
nondecreasing on $[k_0,\iy)$ for some \mbox{$k_0\in\N$}, $\sum
_{k=1}^\infty\Xi_k|\phi''(k)|<\infty$ and\vspace*{1pt} $\limsup_{n\in\N} \Xi
_n|\phi'(n)|=0$ then, for all
$\ell\in\{A,B\}$, $(\Phi_{n,\xi}^\ell)_{n\in\N}$
converges to a finite real
value as $n$ goes to infinity.

In particular, under assumptions \textup{(S)} and \textup{(E2)}, for all
$\ell\!\in\!\{A,B\}$, $(\Phi_{n,\gamma}^\ell)_{n\in\N}$ and\vspace*{2pt}
$(\Phi_{n,\gamma/\G}^\ell)_{n\in\N}$
$[$where $\gamma=(\gamma_n)_{n\in\N}$ and $\gamma/\G=(\gamma
_n/\G_n)_{n\in\N}]$ converge to a finite real
value as $n$ goes to infinity.
\end{Lemma}
\begin{pf}
For all $m$, $n$ $\geq k_0$ with $n\geq m$, Lemma \ref
{Abel1} implies
\[
|\Phi_{n,\xi}^\ell-\Phi_{m,\xi}^\ell|\leq\be_m\Biggl(\sum_{k=m+1}^{n}
\xi_k\phi'(k)+2\xi_m\phi(m)\Biggr).
\]
But
\begin{eqnarray*}
\sum_{k=m+1}^{n}\xi_k\phi'(k)&=&\sum_{k=m+1}^{n}
(\Xi_k-\Xi_{k-1})\phi'(k)
=\sum_{k=m+1}^{n}
\Xi_k\phi'(k)-\sum_{k=m}^{n-1}
\Xi_{k}\phi'(k+1)\\
&=&\sum_{k=m}^{n-1}\Xi_k\bigl(\phi'(k)-\phi'(k+1)\bigr)-\Xi
_m\phi'(m)+\Xi_{n}\phi'(n)\\
&=&-\sum_{k=m}^{n-1}\Xi_k\phi''(k)-\Xi_m\phi'(m)+\Xi_{n}\phi'(n).
\end{eqnarray*}

Let us now prove the convergence of $(\Phi_{n,\gamma}^\ell)_{n\in\N
}$ under assumptions (S)--(E2). Then $\Gamma_n=O(\log n)$ by
Lemma \ref{ubound} and $ \phi'(n)=o(\frac{1}{\log
n})$ (see Remark \ref{ubdphip}) so that $\G_n\phi'(n)\des
_{n\to\iy}0$. Now, there exist $\lambda$, $\mu$ $\in(0,1)$ such that
\begin{eqnarray*}
|\phi''(k)|&=&\bigl|\bigl(\phi(k+1)-\phi(k)\bigr)-\bigl(\phi(k)-\phi(k-1)\bigr)\bigr|\\
&=&
\biggl|\frac{d\phi}{dx}(k+\mu)-\frac{d\phi}{dx}(k-\lambda)\biggr|\\
&\leq&2\sup_{x\in[k-1,k+1]}\biggl|\biggl(\frac{d^2\phi}{dx^2}
\biggr)\biggr|
\end{eqnarray*}
and
\begin{eqnarray*}
\biggl(\frac{d^2\phi}{dx^2}\biggr)(x)&=&\frac{1+\varepsilon
}{(x+2)(\log(x+2))^{2+\varepsilon}}\\
&&{}\times\biggl[-2+\frac{x}{x+2}
\biggl(1+\frac{2+\varepsilon}{\log(x+2)}\biggr)\biggr]\\
&=&O\biggl(\frac{1}{x(\log(x+2))^{2+\e}}\biggr), \qquad  x\in\R
^+\setminus\{0\},
\end{eqnarray*}
so that $\sum\G_k|\phi''(k)|<\iy$ and the assumptions of the first
statement are fulfilled.
The convergence of $(\Phi_{n,\gamma/\G}^\ell)_{n\in\N}$ follows
similarly, since $\gamma_n/\G_n=O(\gamma_n)$.
\end{pf}

%s3.2 ###
\subsection{\texorpdfstring{Proof of Lemma \protect\ref{1sumest}}{Proof of Lemma 6}}
\label{sec:1sumest}
Recall that $\Psi_\iy=0$ [see the first paragraph after the definition of
$(\Psi_n)_{n\in\N}$, Section \ref{sec:sketch}]. Hence, using Lemma
\ref{Abel1},
%
%e13 ###
\begin{eqnarray}\label{psi_maj}
|\Psi_{n}|&=&\Biggl|\sum_{k=n+1}^\infty\frac{\gamma_k}{S_{k-1}}\bigl(\et
_{A,k}-\et_{B,k}-(\tet_A-\tet_B)\bigr)\Biggr|\nonumber\\
&\leq& 2\beta_n\Biggl\{\sum_{k=n+1}^\infty\frac{\gamma
_k}{S_{k-1}}\phi'(k)+2\frac{\gamma_n}{S_{n-1}}\phi(n)\Biggr\}
\\
&\leq& 2\beta_n\Biggl\{\phi'(n)\sum_{k=n+1}^\infty\frac{\gamma
_k}{S_{k-1}}+2\frac{\gamma_n}{S_{n-1}}\phi(n)\Biggr\},\nonumber
\end{eqnarray}
where we use the concavity of $\phi$ in the last inequality.

Now
\[
\sum_{k=n+1}^\iy\frac{\gamma_k}{S_k}=\sum_{k=n+1}^\iy\frac{\De_k}{S_k^2}
=\sum_{k=n+1}^\iy\frac{S_k-S_{k-1}}{S_k^2}\leq\frac{1}{S_n},
\]
so that inequality (\ref{psi_maj}) implies the result.

%s3.3 ###
\subsection{\texorpdfstring{Proof of Lemma \protect\ref{2sumest}}{Proof of Lemma 7}}
\label{sec:2sumest}
Note that
%
%e14 ###
\begin{eqnarray}\label{eq:develop}
\wedge_{n}-\wedge_{m}&=&
\sum_{k=m+1}^n S_{k-1}f(X_{k-1})\frac{\gamma_k}{S_{k-1}}\bigl(\eta
_{A,k}-\eta_{B,k}-(\theta_A-\theta_B)\bigr)\nonumber\\
&=&\sum_{k=m+1}^n
S_{k-1}f(X_{k-1})(\Psi_{k-1}-\Psi_k)\nonumber\\[-8pt]\\[-8pt]
&=&\sum_{k=m+1}^n \Psi_k \bigl(S_kf(X_k)-S_{k-1}f(X_{k-1})\bigr)\nonumber\\
&&{} +\Psi_mS_mf(X_m)-\Psi_nS_nf(X_n).\nonumber
\end{eqnarray}

Recall that $(S_kf(X_k))_{k\in\N}$ is a nondecreasing sequence (see
last paragraph before the statements of Lemmas \ref{1sumest} and \ref
{2sumest}) so that (\ref{eq:develop}) implies, together with
Lem\-ma~\ref{1sumest}, that, for all $n\geq m\geq k_0$,
\begin{eqnarray*}
|{\wedge_{n}-\wedge_{m}}|
&\leq& R'_{m}\Biggl[
\sum_{k=m+1}^n \frac{S_kf(X_k)-S_{k-1}f(X_{k-1})}{S_{k}}
+f(X_m)+f(X_n)
\Biggr]\\
&=& R'_{m}\Biggl[\sum_{k=m+1}^n [f(X_k)-f(X_{k-1})+\gamma_k
f(X_{k-1})]+f(X_m)+f(X_n)\Biggr]\\
&=& R'_{m}\Biggl[
\sum_{k=m+1}^n \gamma_kf(X_{k-1})+2f(X_n)\Biggr].
\end{eqnarray*}

%
%s4 ###
\section{\texorpdfstring{Proof of Theorem \protect\ref{thm:rightconv}}{Proof of Theorem 3}}
%s4.1 ###
\subsection{\texorpdfstring{Brake phenomenon bound: Proof of Proposition \protect\ref{pa=0}}
{Brake phenomenon bound: Proof of Proposition 8}}
\label{sec:pa=0}

Assume that (S) and (E2) hold. Let
\[
\A:=\biggl\{\limsup_{n\to\iy}\frac{Y_n^B}{\log S_n^B}<\iy\biggr\}
\cap\Bigl\{\lim_{n\to\iy}X_n=0\Bigr\}.
\]

In order\vspace*{1pt} to prove Proposition \ref{pa=0}, that is, that $\Pb(\A)=0$,
we first upper bound $S_n^B$ in Lemma \ref{esta}. Then we show that
$Y_n^B\des_{n\to\iy}\iy$ a.s. on $\A$ in Lemma~\ref{ynbiy} so
that, for every $\lambda>0$, $X_n>\lambda/S_n^B$ for large $n\in\N$.
Both lemmas are shown in Section \ref{sec:1pa=0}; we finally
conclude in Section \ref{brakeconc} that $\A$ almost surely does not occur.

%s4.1.1 ###
\subsubsection{Brake phenomenon: Preliminary estimates}
\label{sec:1pa=0}

\begin{Lemma}
\label{esta} Under assumptions \textup{(S)--(E2)}, there exists $L>0$
such that, for all $n\in\N$, $S_n^B\leq Le^{\tet_B\G_n}$ a.s.
\end{Lemma}
\begin{pf}
Recall that (S) implies $\sum\gamma_n^2<\iy$ (see Preliminary
Remark \ref{prrem1}, Section \ref{detest}, or Lemma \ref{ubound}), so that there exists $K>0$ such that,
for all $n\in\N$,
\[
S_n^B\leq K\exp\Biggl(\sum_{k=1}^n\gamma_k\one_{\{I_k>X_k, \eta
_{B,k}=1\}}\Biggr)\qquad\mbox{a.s. }
\]

Now observe that
%
%e15 ###
\begin{eqnarray}\label{decompsum}
\sum_{k=1}^n\gamma_k\one_{\{I_k>X_k,\eta_{B,k}=1\}}
&=&\tet_B\G_n+\sum_{k=1}^n\gamma_k(\eta_{B,k}-\tet_B)
-\sum_{k=1}^n\gamma_k \eta_{B,k}\one_{\{I_k\leq
X_k\}}\nonumber\hspace*{-30pt}\\[-8pt]\\[-8pt]
&=&\theta_B\G_n+\Phi^B_{n,\gamma}-\sum_{k=1}^n \gamma_k\eta
_{B,k}\one_{\{I_k\leq X_k\}},\nonumber
\end{eqnarray}
which enables us to conclude since $\Phi_{n,\gamma}^B$ converges to a
finite value by Lem\-ma~\ref{Abel2}.
\end{pf}
\begin{Lemma}\label{ynbiy}
Under assumptions \textup{(S)--(E2)}, $\A\subseteq\{\limsup_{n\rightarrow
\infty} Y_n^B=\infty\}$, \mbox{$\mathbb P_x$-a.s.}
\end{Lemma}
\begin{pf}
There exist $L$, $L'$ $>0$ such that, for all $n\in\N$,
%
%e16 ###
\begin{equation}
\label{majDeltaB}
\frac{\gamma_{n+1}S_{n}^B}{\G_{n+1}}\leq\frac{\gamma
_{n}S_{n}^B}{\G_{n}}\leq L'e^{-\tet_B\G_n} S_n^B\leq LL',
\end{equation}
where we use (S2) in the first inequality and Lemma \ref{esta}
in the last one.

Now
\begin{eqnarray*}
\Bigl\{\limsup_{n\rightarrow\infty} Y_n^B=\infty\Bigr\}&=&
\Biggl\{\sum_{k=1}^\infty(Y^B_{k+1}-Y_k^B)=\infty\Biggr\}\\
&\supseteq&\Biggl\{\sum_{k=1}^\infty\frac{Y^B_{k+1}-Y_k^B}{\G
_k}=\infty\Biggr\}\\
&=&\Biggl\{\sum_{k=1}^\infty\frac{\Delta_{k+1}^B(1-X_k)}{\G_k}\one
_{\{U_{k+1}=A\}}\eta_{A,k+1}=\infty\Biggr\}\\
&\supseteq&\AAA\cap\Biggl\{\sum_{k=1}^\infty\frac{\gamma
_{k}S_{k-1}^B}{\G_{k}}\one_{\{U_{k}=A\}}\eta_{A,k}=\infty\Biggr\}\\
&=&\AAA\cap\Biggl\{\sum_{k=1}^\infty\frac{\gamma
_{k}S_{k-1}^BX_{k-1}}{\G_{k}}\eta_{A,k}=\infty\Biggr\}\\
&\supseteq&\AAA\cap\Biggl\{\sum_{k=1}^\infty\frac{\gamma_{k}}{\G
_{k}}\eta_{A,k}=\infty\Biggr\}.
\end{eqnarray*}
We use $X_n\des_{n\rightarrow\infty} 0$ a.s. on $\A$ (and
$\gamma_n\to0$) in the second inclusion, whereas, in the third equality,
we apply conditional
Borel--Cantelli lemma (see, e.g.,~\cite{DaDu}, Theorem~2.7.33),\vadjust{\goodbreak}
which claims, given a filtration $\Ff=(\F
_n)_{n\in\N}$ and an $\Ff$-adapted bounded real sequence $(\xi
_n)_{n\geq0}$ (i.e., $\exists M>0$ s.t. $\xi_n\leq M$ a.s.), that
\[
\biggl\{\sum_{n\in\N}\xi_{n}=\iy\biggr\}=\biggl\{\sum_{n\in\N
}\Es(\xi_n|\F_{n-1})=\iy\biggr\}.
\]
Here $\xi_n:=\gamma_{n}S_{n-1}^B\one_{\{U_n=A\}}\eta_{A,n}/\G_n$ is
bounded, using (\ref{majDeltaB}). The last inclusion makes use of
$S_n^BX_n\geq x$ for all $n\in\N$.

Now $\sum\gamma_{k}\eta_{A,k}/\G_{k}=\iy$ a.s. on $\A$, since, on
one hand,
\[
\sum_{k=1}^\infty\frac{\gamma_{k}}{\Gamma_{k}}\geq\sum
_{k=1}^\infty\frac{\Gamma_{k+1}-\Gamma_{k}}{\Gamma_{k}}\geq\int
_{\Gamma_1}^\infty\frac{dx}{x}
\]
and, on the other hand,
\[
\Phi_{n,\gamma/\G}^A:=\sum_{k=1}^n\frac{\gamma_k}{\G_k}(\eta
_{A,k}-\tet_A)
\]
converges (deterministically) to a finite value by Lemma \ref{Abel2}.
\end{pf}

%s4.1.2 ###
\subsubsection{\texorpdfstring{Proof of Proposition \protect\ref{pa=0}}{Proof of Proposition 8}}
\label{brakeconc}

We assume that on the contrary $\mathbb P(\A)>0$ and reach a
contradiction by proving that $\limsup_{n\rightarrow\infty}
Y_n^B/\log(S_n^B)=\infty$ a.s. on $\A$.
Note that
\[
Y_n^B=\sum_{k=0}^{n-1}\Delta_{k+1}^B \one_{\{I_{k+1}\leq X_k\}}\eta
_{A,k+1}(1-X_k)+x
\]
and let, for all $\lambda>0$,
\begin{eqnarray*}
Z_n^{B,\lambda}&:=&\sum_{k=0}^{n-1}\gamma_{k+1}S_{k}^B \one_{\{
I_{k+1}\leq\lambda/S_{k}^B\}}\eta_{A,k+1},\\
{\tilde Z}_n^{B,\lambda}&:=&\sum_{k=0}^{n-1}\gamma_{k+1}S_{k}^B \min
\biggl(1,\frac{\lambda}{S_{k}^B}\biggr)\eta_{A,k+1}.
\end{eqnarray*}
Almost surely on $\A$, $\limsup_{n\rightarrow\infty} Y_n^B=\infty$
by Lemma \ref{ynbiy} and $\lim_{n\rightarrow\infty} X_n=\lim
_{n\rightarrow\infty}\gamma_n=0$, so that, for all $\lambda>0$
\[
\limsup_{n\rightarrow\infty} \frac{Y_n^B}{\log(S_n^B)}\geq\limsup
_{n\rightarrow\infty} \frac{Z_n^{B,\lambda}}{\log(S_n^B)}  \qquad
\mbox{a.s.}
\]
Fix $\lambda>0$. To show that the right-hand side of this last inequality is infinite
a.s. on~$\A$, we aim to estimate $\Es(Z_n^{B,\lambda})=\Es({\tilde
Z}_n^{B,\lambda})$ and to upper bound $\Es((Z_n^{B,\lambda}-{\tilde
Z}_n^{B,\lambda})^2)$. In order to yield the latter we first observe
that there exists $M>0$ such that, for all $k\in\N$, $\gamma
_{k+1}S_k^B\leq\De_k^B\leq M\G_k$, by
inequali\-ty~(\ref{majDeltaB}).\looseness=1

Now
%
%e17 ###
\begin{eqnarray}\label{ubvar}
&&\Es\bigl((Z_n^{B,\lambda}-{\tilde Z}_n^{B,\lambda})^2\nonumber\bigr)\\
&&\qquad=\mathbb E\Biggl( \sum_{k=0}^{n-1}(\gamma_{k+1}S_k^B)^2
\min\biggl(1,\frac{\lambda}{S_{k}^B}\biggr)
\biggl(1- \min\biggl(1,\frac{\lambda}{S_{k}^B}\biggr)\biggr)\eta
_{A,k+1}\Biggr)\\
&&\qquad\leq M\G_n\mathbb E\Biggl( \sum_{k=0}^{n-1} \gamma_{k+1}S_k^B\min
\biggl(1,\frac{\lambda}{S_{k}^B}\biggr)\eta_{A,k+1}\Biggr)
=M \G_n \mathbb E(Z_{n}^{B,\lambda}).\nonumber
\end{eqnarray}

On the other hand, for all $M>0$ and $\e>0$,
\begin{eqnarray*}
\Es(Z^{B,\lambda}_{n})&=&\Es\Biggl( \sum_{k=0}^{n-1}\gamma
_{k+1}S_k^B \min\biggl(1,\frac{\lambda}{S_{k}^B}\biggr)\eta
_{A,k+1} \Biggr)\\
&\geq&\lambda(1-\e)\Pb(\A)\sum_{k=k_0(\e,\lambda)}^{n-1}
\gamma_{k+1}\eta_{A,k+1},
\end{eqnarray*}
where we use that $S_n^B=Y_n^B/X_n\to\iy$ a.s. on $\A$, $k_0(\e
,\lambda)$ being a constant depending on $\e$ and $\lambda$. Now
$\Phi_{n,\gamma}^A=\sum_{k=0}^{n-1}
\gamma_{k+1}\eta_{A,k+1}-\G_n\tet_A$ converges by Lemma \ref
{Abel2}, so that we obtain
\[
\lambda\tet_A\geq\limsup_{n\to\iy}\frac{\Es(Z_{n}^{B,\lambda
})}{\G_n}\geq\liminf_{n\to\iy}\frac{\Es(Z_{n}^{B,\lambda})}{\G
_n}\geq\lambda\Pb(\A)\tet_A.
\]
Fix $\rho\in(0,1)$ and let
\[
B_{n,\lambda}:=\{|Z_{n}^{B,\lambda}-{\tilde
Z}_n^{B,\lambda}|\leq\rho\mathbb E(Z_{n}^{B,\lambda})\}.
\]
By (\ref{ubvar}) and Chebyshev's inequality,
\[
\Pb(B_{n,\lambda}^c)\leq\frac{M\G_n}{\rho^2 \mathbb E
(Z_{n}^{B,\lambda})}.
\]
Therefore, for all $\lambda>0$, if we let $\C_\lambda:=\A\cap
\limsup_{n\to\iy} B_{n,\lambda}$,
\[
\Pb(\C_\lambda)\geq\limsup_{n\to\iy}\Pb(\A\cap B_{n,\lambda
})\geq\Pb(\A)-\frac{M}{\lambda\rho^2\tet_A\Pb(\A)}>0,
\]
if we choose $\lambda$ such that
$\lambda>M\tet_A^{-1}(\rho\Pb(\A))^{-2}$.

Now, almost surely on $\C_\lambda\subseteq\A$, ${\tilde
Z}_n^{B,\lambda}/\G_n\des_{n\to\iy}\lambda\tet_A$ (since
$S_n^B\des_{n\to\iy}\iy$; see above), so that
\[
\limsup_{n\to\iy}\frac{Y_n^B}{\log S_n^B}\geq\frac{\lambda(1-\rho)
\tet_A}{\tet_B}
\]
using that $\limsup_{n\to\iy}\log S_n^B/\G_n\leq\tet_B$ by Lemma
\ref{esta}.

Therefore,
\begin{eqnarray*}
\Pb\biggl(\biggl\{\limsup_{n\to\iy}\frac{Y_n^B}{\log S_n^B}=\iy
\biggr\}\cap\A\biggr)
&\geq&\Pb\Bigl(\limsup_{\lambda\in\N, \lambda\to\iy}\C_\lambda
\Bigr)\\
&\geq&\limsup_{\lambda\in\N, \lambda\to\iy}\Pb(\C_\lambda)\geq
\Pb(\A),
\end{eqnarray*}
which enables us to conclude.

%s4.2 ###
\subsection{\texorpdfstring{Conclusion of the proof of Theorem \protect\ref{thm:rightconv}}{Conclusion of the proof of Theorem 3}}
\label{sec:conc}

Let, for all $n\geq0$, $T_n^B:=e^{\tet_B \G_n}$. It follows from
Proposition \ref{pa=0}
that
\[
\limsup_{n\rightarrow\infty} \frac{X_n}{\log T_n^B/T_n^B}=\infty
\qquad\mbox{a.s. on }X_\iy=0
\]
using that $\limsup_{n\rightarrow\infty}S_n^B/T_n^B<\iy$ by Lemma
\ref{esta}.

Given $l\in\N$, let us estimate $\Pb(X_\infty=0|\F_l)$. Using
identity (\ref{evolm}) and the assumption $\tet_A>\tet_B$, there
exists $n_0\in\N$ deterministic such that, for all $n\geq m\geq n_0$,
\begin{eqnarray*}
X_n-X_m&=&M_n-M_m+\bigl(\tet_A-\tet_B+\Box(R'_m)\bigr)\sum_{k=m+1}^n \gamma
_kf(X_{k-1})\\
&&{} +2\Box(R'_m) f(X_n)\\
&\geq& M_n-M_m-X_n,
\end{eqnarray*}
so that
%
%e18 ###
\begin{equation}
\label{evol0}
2X_n\geq X_m+M_n-M_m.
\end{equation}

Let $(N_n)_{n\geq l}$ be the $(\F_n)_{n\geq l}$ adapted martingale
given by
\[
N_n:=\sum_{i=l+1}^n \gamma_i\one_{\{X_{i-1}\leq
X_l\}}\varepsilon_i,\qquad
N_l:=0;
\]
recall that $(\e_i)_{i\in\N}$ was defined before the statement of
Proposition \ref{evol}.

Let $n_0$ be sufficiently large, so that $\gamma_{n_0}\leq1/2$; then,
for all $n\geq n_0$, $X_{n+1}>X_n/2$. Thus, for all
$n\geq l\geq n_0$, inequality (\ref{evol0}) implies
%
%e19 ###
\begin{equation}
\label{evol1}
2X_n\geq X_m+N_n-N_m\geq X_l/2+N_n-N_m,
\end{equation}
where $m:=\max\{l\leq i\leq n\dvtx X_i>X_l/2\}$; indeed, if
$m<n$ then, for all $m\leq k\leq n-1$, $X_{k+1}\leq X_l/2$, hence, $X_k
\leq X_l$; (\ref{evol1}) also trivially holds in the case $n=m$. Hence,
if $x^-:=\max(-x,0)$ denotes the negative part of $x$, then
\[
(2X_\infty-X_l/2)^-\leq\sup_{m,n\geq l}|N_n-N_m|\leq2\sup_{n\geq l}|N_n-N_l|.
\]
Therefore, by Chebyshev's inequality,
%
%e20 ###
\begin{eqnarray}\label{majdeltalim1}
\Pb(X_\infty=0|\F_l)&\leq&\frac{4\Eb[[(2X_\infty
-X_l/2)^-]^2|\F_l]}{X_l^2}\nonumber\\[-8pt]\\[-8pt]
&\leq&16 \frac{\Eb[\sup_{n\geq l}(N_n-N_l)^2|\F_l]}{X_l^2}.\nonumber
\end{eqnarray}
Now observe that, for all $k\in\N$, $\Eb(\varepsilon^2_{k+1}|\F
_k)\leq f(X_k)\leq X_k$, so that Doob's inequality implies
%
%e21 ###
\begin{eqnarray}\label{majdeltalim2}
\Eb\Bigl[\sup_{n\geq l}(N_n-N_l)^2\big|\F_l\Bigr]
&\leq&4\Eb\Biggl(\sum_{n=l+1}^\infty\gamma_n^2 \one_{\{X_{n-1}\leq X_l\}
}f(X_{n-1})\Biggr)\nonumber\\[-8pt]\\[-8pt]
&\leq&4X_l\sum_{n=l+1}^\infty\gamma_n^2.\nonumber
\end{eqnarray}
Let us upper bound $\sum_{i=n+1}^\infty\gamma_i^2$ in terms of
$T_n$. For sufficiently large $k\in\N$,
\[
T_{k+1}^B-T_k^B=e^{\tet_B\G_{k+1}}(1-e^{-\tet_B \gamma_{k+1}})\geq
\frac{T_{k+1}^B\tet_B\gamma_{k+1}}{2}
\]
and, on the other hand, by assumption (S),
\[
\gamma_k\leq C\Gamma_k e^{-\tet_B \Gamma_{k}}=\frac{C\log
(T_k^B)}{\tet_B T_k^B}.
\]
Hence, if $l\in\N$ was assumed sufficiently large,
%
%e22 ###
\begin{equation}\label{newmajsigmsq}\qquad
\sum_{n=l+1}^\infty\gamma_n^2\leq C\sum_{n=l+1}^\infty
(T_n^B-T_{n-1}^B) \frac{\log T^B_n}{(T^B_n)^2}
\leq C\int_{T_l^B}^\infty\frac{\log t}{t^2}\,dt\leq 2C\frac{\log T^B_l}{T_l^B}.
\end{equation}
In summary, it follows from identities (\ref{majdeltalim1})--(\ref
{newmajsigmsq}) that
\[
\Pb(X_\infty=0|\F_l)\leq C\frac{\log T^B_l}{X_lT_l^B}.
\]
Now the bounded martingale $\Pb(X_\infty=0|\F_l)$ converges, as $l$
goes to infinity,~to
\[
\one_{\{X_\infty=0\}}\leq C\liminf_{l\rightarrow\infty} \frac{\log
T^B_l}{X_lT_l^B}=0\qquad\mbox{a.s.}
\]
so that $\Pb(X_\infty=0)=0$.

\section*{Acknowledgments}
We are grateful to the referee for very helpful comments.

%suskaldyti doi

% imsref loaded by lrinkeviciute, 2011-04-07 08:58:49
% imsref loaded by lrinkeviciute, 2011-04-07 09:01:39

%
\printaddresses

\end{document}